\documentclass[12pt]{amsart}
\usepackage{amsmath,amsthm,amssymb,amscd}
\usepackage{graphicx, psfrag}

\newdimen\tableauside\tableauside=1.0ex
\newdimen\tableaurule\tableaurule=0.4pt
\newdimen\tableaustep
\def\phantomhrule#1{\hbox{\vbox to0pt{\hrule height\tableaurule width#1\vss}}}
\def\phantomvrule#1{\vbox{\hbox to0pt{\vrule width\tableaurule height#1\hss}}}
\def\sqr{\vbox{%
  \phantomhrule\tableaustep
  \hbox{\phantomvrule\tableaustep\kern\tableaustep\phantomvrule\tableaustep}%
  \hbox{\vbox{\phantomhrule\tableauside}\kern-\tableaurule}}}
\def\squares#1{\hbox{\count0=#1\noindent\loop\sqr
  \advance\count0 by-1 \ifnum\count0>0\repeat}}
\def\tableau#1{\vcenter{\offinterlineskip
  \tableaustep=\tableauside\advance\tableaustep by-\tableaurule
  \kern\normallineskip\hbox
    {\kern\normallineskip\vbox
      {\gettableau#1 0 }%
     \kern\normallineskip\kern\tableaurule}%
  \kern\normallineskip\kern\tableaurule}}
\def\gettableau#1 {\ifnum#1=0\let\next=\null\else
  \squares{#1}\let\next=\gettableau\fi\next}

 \numberwithin{equation}{section}

 \newcommand{\R}{\mathbb{R}}
 \newcommand{\C}{\mathbb{C}}
 \newcommand{\Z}{\mathbb{Z}}
 \newcommand{\CC}{\mathcal{C}}
 \newcommand{\CH}{\mathcal{H}}
 \newcommand{\CL}{\mathcal{L}}
 \newcommand{\g}{\mathfrak{g}}
 \newcommand{\Ug}{{U_q(\g)}}
 \newcommand{\SU}{\mathrm{SU}}
 \newcommand{\slN}{\mathrm{sl}_N}
 \newcommand{\UslN}{{U_q(\slN)}}
 \newcommand{\cR}{\check{R}}
 \newcommand{\otn}{{\otimes n}}
 \newcommand{\otr}{{\otimes r}}
 \newcommand{\V}[1]{{V^{\otimes #1}}}
 \newcommand{\Vn}{{V^{\otimes n}}}
 \newcommand{\x}{\mathbf{x}}

 \DeclareMathOperator{\id}{id}
 \DeclareMathOperator{\lk}{lk}
 \DeclareMathOperator{\tr}{tr}
 \DeclareMathOperator{\End}{End}
 \DeclareMathOperator{\Hom}{Hom}
 \DeclareMathOperator{\Ker}{Ker}

 \newtheorem{thm}{Theorem}[section]
 \newtheorem{lem}[thm]{Lemma}
 \newtheorem{prop}[thm]{Proposition}
 \newtheorem{cor}[thm]{Corollary}
 
 \newtheorem{conj}[thm]{Conjecture}
 \theoremstyle{definition}
 \newtheorem{defn}[thm]{Definition}
 \newtheorem{exam}[thm]{Example}
 \newtheorem{rem}[thm]{Remark}
 \theoremstyle{remark}

\begin{document}

\title[Hecke algebra and colored HOMFLY polynomial]{On the Hecke algebras and \\ the colored HOMFLY polynomial}
\author[Xiao-Song Lin]{Xiao-Song Lin}
\address{Department of Mathematics, University of California,
Riverside, CA 92521, USA}
\thanks{The first author is supported in part by NSF grants DMS-0404511}
\email{xl@math.ucr.edu}
\author[Hao Zheng]{Hao Zheng}
\address{Department of Mathematics, Zhongshan University, Guangzhou 510275, China}
\thanks{}
\email{zhenghao@mail.sysu.edu.cn}
\date{}

\begin{abstract}
The colored HOMFLY polynomial is the quantum invariant of oriented
links in $S^3$ associated with irreducible representations of the
quantum group $U_q(\mathrm{sl}_N)$. In this paper, using an
approach to calculate quantum invariants of links via
cabling-projection rule, we derive a formula for the colored
HOMFLY polynomial in terms of the characters of the Hecke algebras
and Schur polynomials. The technique leads to a fairly simple
formula for the colored HOMFLY polynomial of torus links. This
formula allows us to test the Labastida-Mari\~no-Vafa conjecture,
which reveals a deep relationship between Chern-Simons gauge
theory and string theory, on torus links.
\end{abstract}
\maketitle

\section{Introduction}

In the abstract of his seminal paper \cite{Jones}, V. Jones wrote:
``By studying representations of the braid group satisfying a
certain quadratic relation we obtain a polynomial invariant in two
variables for oriented links. \dots The two-variable polynomial
was first discovered by Freyd-Yetter, Lickorish-Millet, Ocneanu,
Hoste, and Przytycki-Traczyk.'' This two variable link polynomial
$P_\CL(t,\nu)$, commonly referred to as the HOMFLY polynomial for
an oriented link $\CL$ in $S^3$, is characterized by the following
crossing changing formula:
\begin{eqnarray}
  && P_\text{unknot}(t,\nu) = 1, \\
  && \nu^{-1/2} P_{\CL_+}(t,\nu) - \nu^{1/2} P_{\CL_-}(t,\nu)
  = (t^{-1/2}-t^{1/2}) P_{\CL_0}(t,\nu).
\end{eqnarray}

Since then, this two variable link polynomial has been generalized
to the quantum invariant associated with irreducible
representations of the quantum group $\UslN$, with the variables
$t^{1/2}=q^{-1}$ and $\nu^{1/2}=q^{-N}$. We will refer to this
generalization as the colored HOMFLY polynomial.

Despite the fact that the theory of quantum invariants of links is
by now well developed, the computation of colored HOMFLY
polynomial is still extremely challenging. Besides the trivial
links, a general formula seems to exist in the mathematics
literature only for the Hopf link \cite{Morton-Lukac}. In the
physics literature, Witten's Chern-Simons path integral with the
gauge group $\SU_N$ \cite{Witten} offers an intrinsic but not
rigorous definition of the colored HOMFLY polynomial. There is a
conjectured relationship between the $1/N$ expansion of
Chern-Simons theory and the Gromov-Witten invariants of certain
non-compact Calabi-Yau 3-folds. See
\cite{Gopakumar-Vafa}\cite{Ooguri-Vafa} for example. Motivated by
this conjectured relationship, Labastida, Mari\~no and Vafa
proposed a precise conjecture about the structure of their
reformulation of the colored HOMFLY polynomial
\cite{Labastida-Marino2}\cite{LMV}. See Section 5. A formula of
the colored HOMFLY polynomial for torus knots is given in
\cite{Labastida-Marino}, which was used to test the
Labastida-Mari\~no-Vafa conjecture on torus knots.

In this paper, using an approach to calculate quantum invariants
of links via cabling-projection rule, we derive a formula for the
colored HOMFLY polynomial in terms of the characters of the Hecke
algebras and Schur polynomials. See Theorem 4.3. An important
feature of this formula is that the character of the Hecke algebra
is free of the variable $\nu$ and the Schur polynomial is
independent of the link $\CL$. We think that this separation of
the variable $\nu$ and the link $\CL$ might be important for a
possible proof of the Labastida-Mari\~no-Vafa conjecture.

Our technique leads to a fairly simple formula for the colored
HOMFLY polynomial of torus links. See Theorem 5.1. Using our formula,
the
Labastida-Mari\~no-Vafa conjecture can be test on 
several
infinite families of torus links. 
Our calculation also suggests 
a new
structure of the reformulated colored HOMFLY polynomial of torus
links: it is equivalent to a family of polynomials in $\Z[t^{\pm1}]$
invariant under the transformation $t\rightarrow t^{-1}$. 
See Conjecture 6.2 and the examples following it.

\medskip

\noindent{\bf Acknowledgments.} The authors would like to thank
Professors Kefeng Liu and Jian Zhou for their interest in this
work.

\section{Link invariants from quantum groups}\label{sec:inv}

In this section, we give a brief review of the quantum group
invariants of links. See \cite{Kassel}\cite{RT}\cite{Tur} for
details.

Let $\g$ be a complex simple Lie algebra and let $q$ be a nonzero
complex number which is not a root of unity. Let $\Ug$ denote the
quantum enveloping algebra of $\g$. The ribbon category structure
of the set of finite dimensional complex representations of $\Ug$
provides the following objects.

1. Associated to each pair of $\Ug$-modules $V,W$ is a natural
isomorphism (the braiding) $\cR_{V,W} : V \otimes W \to W \otimes
V$ such that
\begin{equation}\label{eqn:cR}
\begin{split}
  & \cR_{U \otimes V,W} =
  (\cR_{U,W} \otimes \id_V)(\id_U \otimes \cR_{V,W}), \\
  & \cR_{U,V \otimes W} =
  (\id_V \otimes \cR_{U,W})(\cR_{U,V} \otimes \id_W)
\end{split}
\end{equation}
hold for all $\Ug$-modules $U,V,W$. The naturality means
\begin{equation}
  (y \otimes x) \cR_{V,W} = \cR_{V',W'} (x \otimes y)
\end{equation}
for $x \in \Hom_\Ug(V,V')$, $y\in \Hom_\Ug(W,W')$. These
equalities imply the braiding relation
\begin{equation}\label{eqn:brading}
\begin{split}
  (\cR_{V,W}\otimes \id_U) & (\id_V \otimes \cR_{U,W}) (\cR_{U,V} \otimes \id_W) \\
  & = (\id_W \otimes \cR_{U,V}) (\cR_{U,W} \otimes \id_V) (\id_U \otimes \cR_{V,W}).
\end{split}
\end{equation}

2. There exists an element $K_{2\rho} \in \Ug$ (the enhancement of
$\cR$, here $\rho$ means the half-sum of all positive roots of
$\g$) such that
\begin{equation}\label{eqn:k2rho}
  K_{2\rho} (v \otimes w)
  = K_{2\rho}(v) \otimes K_{2\rho}(w)
\end{equation}
for $v \in V$, $w \in W$. Moreover, for every $z \in \End_\Ug(V
\otimes W)$ with $z = \sum_i x_i \otimes y_i$, $x_i \in \End(V)$,
$y_i \in \End(W)$ one has the {\em (partial) quantum trace}
\begin{equation}
  \tr_W(z) = \sum_i \tr(y_iK_{2\rho}) \cdot x_i \in \End_\Ug(V).
\end{equation}

3. Associated to each $\Ug$-module $V$ is a natural isomorphism
(the ribbon structure) $\theta_V : V \to V$ satisfying
\begin{equation}\label{eqn:theta}
  \theta_V^{\pm1} = \tr_V \cR_{V,V}^{\pm1}.
\end{equation}
The naturality means
\begin{equation}
  x \cdot \theta_V = \theta_{V'} \cdot x
\end{equation}
for $x \in \Hom_\Ug(V,V')$.

\bigskip
\begin{center}
    \psfrag{V}[]{$V$}
    \psfrag{W}[]{$W$}
    \psfrag{L}[]{$\CL$}
    \psfrag{b}[]{$\beta$}
    \includegraphics[scale=.8]{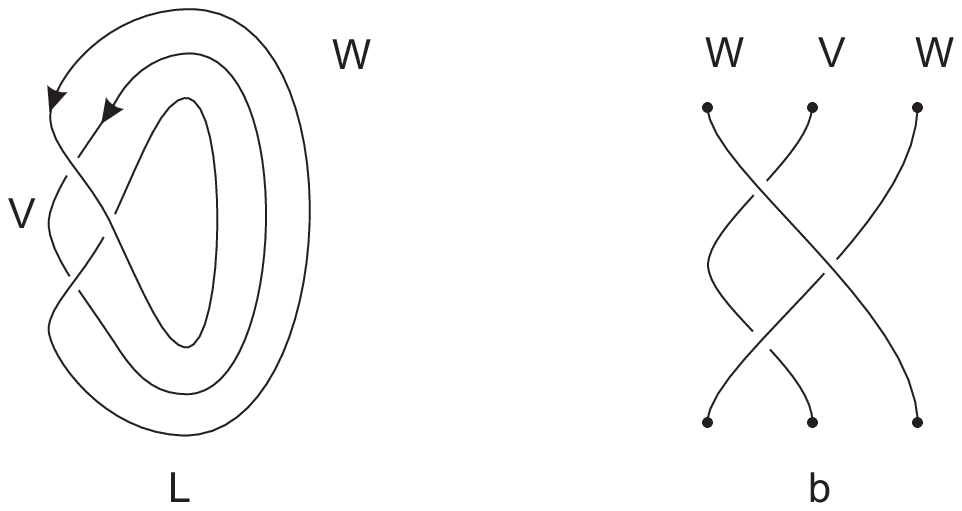}
\end{center}
\bigskip

With these objects, one constructs the quantum group invariants of
links as follows. Let $\CL$ be an oriented link with the
components $\CL_1,\dots,\CL_l$ labeled by the $\Ug$-modules
$V_1,\dots,V_l$, respectively. Choose a closed braid
representative $\hat\beta$ of $\CL$ with $\beta \in B_n$ being an
$n$-strand braid. Assign to each positive (resp. negative)
crossing of $\beta$ an isomorphism $\cR_{V,W}$ (resp.
$\cR_{W,V}^{-1}$) where $V,W$ are the $\Ug$-modules labeling the
two outgoing strands of the crossing.

\bigskip
\begin{center}
    \psfrag{V}[]{$V$}
    \psfrag{W}[]{$W$}
    \psfrag{P}[]{$\cR_{V,W}$}
    \psfrag{N}[]{$\cR_{W,V}^{-1}$}
    \includegraphics[scale=.7]{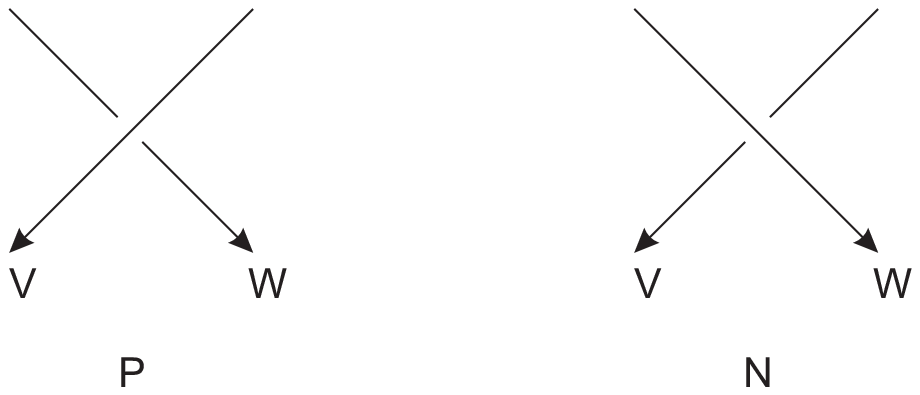}
\end{center}
\bigskip

Then the braid $\beta$ gives rise to an isomorphism
\begin{equation}
  h_{V'_1,\dots,V'_n}(\beta) \in \End_\Ug(V'_1 \otimes \cdots \otimes V'_n),
\end{equation}
where $V'_1,\dots,V'_n$ are the $\Ug$-modules labeling the strands
of $\beta$, and the quantum trace
\begin{equation}
  \tr_{V'_1 \otimes \cdots \otimes V'_n} h_{V'_1,\dots,V'_n}(\beta)
\end{equation}
defines a framing dependent link invariant of $\CL$.

\begin{exam}
The link shown in above figure has two components, labeled by $W$
and $V$ respectively. It is the closure of $\beta =
\sigma_1^{-1}\sigma_2^{-1}\sigma_1 \in B_3$, which gives rise to
an isomorphism
\begin{equation}
  h_{W,V,W}(\beta)
  = (\cR_{W,V}^{-1}\otimes\id_W)
  (\id_V\otimes\cR_{W,W}^{-1}) (\cR_{W,V}\otimes\id_W).
\end{equation}
Thus the link invariant is
\begin{equation}
  \tr_{W \otimes V \otimes W}
  (\cR_{W,V}^{-1}\otimes\id_W)
  (\id_V\otimes\cR_{W,W}^{-1}) (\cR_{W,V}\otimes\id_W).
\end{equation}
\end{exam}

To eliminate the framing dependency, one should require the
modules $V_1,\dots,V_l$ be irreducible, hence the isomorphisms
$\theta_{V_1},\dots,\theta_{V_l}$ are multiples of identity and
may be regarded as scalars. Let $w(\CL_i)$ be the {\em writhe} of
$\CL_i$ in $\beta$, i.e. the number of positive crossings minus
the number of negative crossings. Then the quantity
\begin{equation}\label{eqn:inv_defn}
  I_{\CL;V_1,\dots,V_l} =
  \theta_{V_1}^{-w(\CL_1)} \cdots \theta_{V_l}^{-w(\CL_l)}
  \tr_{V'_1 \otimes \cdots \otimes V'_n} h_{V'_1,\dots,V'_n}(\beta)
\end{equation}
defines a framing independent link invariant.

When the link involved is the unknot, it is easy to see that
\begin{equation}
  I_{\text{unknot};V} = \tr_V \id_V.
\end{equation}
This quantity is regarded as the quantum version of the classical
dimension of $V$, referred to as the {\em quantum dimension} of
$V$ and denoted by $\dim_q V$.

\section{Centralizer algebra and cabling-projection rule}\label{sec:central}

In general, the isomorphism $\cR_{V,W}$ is very complicated when
the dimensions of $V,W$ are larger, so it is not practical to
compute the link invariants from their definition. However, on the
other hand, general representations of a simple Lie algebra $\g$
(thus its quantum deformation $\Ug$) are often realized as
components of tensor products of some simple ones. For example,
irreducible representations of $\UslN$ are always the components
of some tensor products of the fundamental representation.

In this section, we follow this observation and develop a
cabling-projection rule to break down the complexity of general
$\cR$. For this purpose we need the notion of centralizer algebra.

The centralizer algebras of the modules of simple Lie algebras
have played an important role in representation theory. Parts of
their quantum version were studied in
\cite{BW}\cite{LR}\cite{Wen}. In the case of $\UslN$, the
situation is desirable. The centralizer algebras are nothing but
the subalgebras of the Hecke algebras of type $A$.

\medskip

Let $V$ be a $\Ug$-module. The {\em centralizer algebra} of $\Vn$
is defined as
\begin{equation}
  \CC_n(V) = \End_\Ug(\Vn) =
  \{ x \in \End(\Vn) \mid xy = yx, \; \forall y \in \Ug \}.
\end{equation}
It is immediate from definition that $\CC_n(V)$ is a finite
dimensional von Neumann algebra, i.e. the algebra is isomorphic to
a direct sum of matrix algebras. Indeed, if $\Vn$ admits the
irreducible decomposition
\begin{equation}\label{eqn:decom_vn}
  \Vn = \bigoplus_{\lambda \in \Lambda} d_\lambda \cdot V_\lambda,
\end{equation}
by Schur's lemma we have
\begin{equation}\label{eqn:decom_cn}
  \CC_n(V) = \bigoplus_{\lambda \in \Lambda} \CC_\lambda
\end{equation}
where $\CC_\lambda = \End_\Ug(d_\lambda V_\lambda)$ is a full
$d_\lambda \times d_\lambda$ matrix algebra. Since each matrix
algebra admits a unique irreducible representation, via above
decomposition the irreducible representations of $\CC_n(V)$ are
naturally indexed by $\Lambda$.

Let $\zeta^\lambda$ denote the character of the irreducible
representation of $\CC_n(V)$ indexed by $\lambda \in \Lambda$.

\begin{lem}\label{lem:trace}
For every $x \in \CC_n(V)$ we have
\begin{equation}
  \tr_\Vn x
  = \sum_{\lambda \in \Lambda}
  \zeta^\lambda(x) \cdot \dim_q V_\lambda.
\end{equation}
\end{lem}

\begin{proof}
Let $\pi_\lambda$ be the unit of $\CC_\lambda$. Then $\pi_\lambda
x$ is a matrix in $\CC_\lambda$, whose normal trace $\tr
\pi_\lambda x$ is precisely $\zeta^\lambda(x)$. Therefore,
\begin{equation}
  \tr_\Vn x
  = \sum_{\lambda \in \Lambda}
  \tr \pi_\lambda x \cdot \tr_{V_\lambda} \id_{V_\lambda}
  = \sum_{\lambda \in \Lambda}
  \zeta^\lambda(x) \cdot \dim_q V_\lambda.
\end{equation}
\end{proof}
A {\em projection (or idempotent)} of $\CC_n(V)$ is an element $p
\in \CC_n(V)$ satisfying the idempotent equation $p^2 = p$. By
definition, an element $p \in \CC_n(V)$ is a projection if and
only if it is, restricted on each $\CC_\lambda$, diagonalizable
and has the only possible eigenvalues $0$ and $1$. It is clear
that for each projection $p \in \CC_n(V)$,
\begin{equation}\label{eqn:proj_decom}
  p\Vn \cong
  \bigoplus_{\lambda \in \Lambda}
  \zeta^\lambda(p) \cdot V_\lambda.
\end{equation}
A projection $p \in \CC_n(V)$ is called {\em minimal (or
primitive)} if  $p \Vn \cong V_\lambda$ for some $\lambda \in
\Lambda$.

Let $h_V$ be the homomorphism
\begin{equation}
  h_V : \C B_n \to \CC_n(V), \quad
  \sigma_i \mapsto
  \id_\V{(i-1)} \otimes R_{V,V} \otimes \id_\V{(n-i-1)}.
\end{equation}
The following lemma makes it possible to recover general $\cR$,
$\theta$ from specific ones.

\bigskip
\begin{center}
    \psfrag{D}[]{$\Delta_5^2$}
    \psfrag{X}[]{$\chi_{3,4}$}
    \includegraphics[scale=.9]{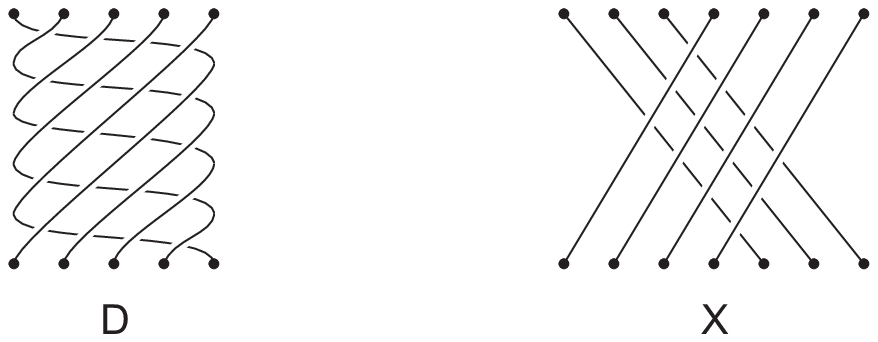}
\end{center}
\bigskip

\begin{lem}\label{lem:theta}
Let $\Delta_n^2 = (\sigma_1\sigma_2\cdots\sigma_{n-1})^n \in B_n$
and $\chi_{n,n'} = \prod_{i=1}^{n'}
(\sigma_{i+n-1}\sigma_{i+n-2}\cdots\sigma_i) \in B_{n+n'}$ denote
the full twist braid and the $(n,n')$-crossing braid,
respectively. Let $p \in \CC_n(V)$, $p' \in \CC_{n'}(V)$ be
projections and $U=p\Vn$, $W=p'\V{n'}$. We have
\begin{eqnarray}
  \label{eqn:tensorR}
  && h_V(\chi_{n,n'}) \cdot (p \otimes p')
  = \cR_{W,U} \oplus 0_{\Ker p \otimes p'}, \\
  && (\theta_V)^\otn \cdot h_V(\Delta_n^2) \cdot p
  = \theta_{U} \oplus 0_{\Ker p}.
\end{eqnarray}
\end{lem}

\begin{proof}
Applying the identities (\ref{eqn:cR}), (\ref{eqn:brading}) and
(\ref{eqn:theta}) inductively, we have
\begin{eqnarray}
  && \cR_{\V{n'},\Vn} = h_V(\chi_{n,n'}), \\
  && \theta_{\Vn} = \tr_\Vn \cR_{\Vn,\Vn}
  = (\theta_V)^\otn \cdot h_V(\Delta_n^2).
\end{eqnarray}
Then from the naturality of $\cR$ and $\theta$ the lemma follows.
\end{proof}

As an easy consequence of (\ref{eqn:tensorR}) and the naturality
of $\cR$, we have the cabling-projection rule

\begin{lem}\label{lem:cable}
Let $\beta \in B_m$ be a braid and $p_i \in \CC_{n_i}(V)$,
$i=1,\dots,m$ be projections such that $p_i=p_j$ whenever the
$i$-th strand of $\beta$ ends at $j$-th point. Moreover, let
$V_i=p_i \V{n_i}$, $n=n_1+\cdots+n_m$ and $\beta^{(n_1,\dots,n_m)}
\in B_n$ be the braid obtained by cabling the $i$-th strand of
$\beta$ to $n_i$ parallel ones. Then
\begin{equation}
  h_V(\beta^{(n_1,\dots,n_m)}) \cdot (p_1 \otimes \cdots \otimes p_m)
  = h_{V_1,\dots,V_m}(\beta) \oplus
  0_{\Ker p_1 \otimes \cdots \otimes p_m},
\end{equation}
thus
\begin{equation}
  \tr_{V_1 \otimes \cdots \otimes V_m} h_{V_1,\dots,V_m}(\beta)
  = \tr_\Vn h_V(\beta^{(n_1,\dots,n_m)}) \cdot
  (p_1 \otimes \cdots \otimes p_m).
\end{equation}
\end{lem}

With above lemmas, one is able to re-express the link invariant
(\ref{eqn:inv_defn}), by choosing a suitable $\Ug$-module $V$, in
terms of much more accessible objects: the characters and
projections of the centralizer algebras $\CC_{n}(V)$ and the
quantum traces of $\Ug$-modules. In the next section, we present a
detailed realization of this approach for the case $\g = \slN$.

\section{Hecke algebras and colored HOMFLY polynomial}\label{sec:homfly}

In the rest part of this paper we will extensively apply the facts
concerning the Hecke algebras, the quantum enveloping algebras
$\UslN$ and symmetric functions. The facts are well known and most
of them can be found, for example, in
\cite{DJ1}\cite{Kassel}\cite{KS}\cite{Mur}\cite{Sag}.

The Hecke algebra $\CH_n(q)$ of type $A_{n-1}$ is the complex
algebra with generators $g_1,g_2,\dots,g_{n-1}$ and relations
\begin{equation}
  \begin{array}{ll}
    g_i g_j = g_j g_i, & |i-j| \geq 2, \\
    g_ig_jg_i = g_jg_ig_j, & |i-j| = 1, \\
    (g_i-q) (g_i+q^{-1}) = 0, & i=1,2,\dots,n-1. \\
  \end{array}
\end{equation}
Note that, when $q=1$, the Hecke algebra $\CH_n(q)$ is nothing new
but the group algebra $\C\Sigma_n$ of the symmetric group. In
fact, if $q$ is nonzero and not root of unity we still have the
isomorphism $\CH_n(q) \cong \C\Sigma_n$ and $\CH_n(q)$ also
canonically decomposes as
\begin{equation}\label{eqn:decom_hn}
  \CH_n(q) = \bigoplus_{\lambda \vdash n} \CH_\lambda(q)
\end{equation}
with each $\CH_\lambda(q)$ being a matrix algebra.

Here we fix several notations of combinatorics. A {\em
composition} $\mu$ of $n$, denoted by $\mu \models n$, is a
sequence of nonnegative integers $(\mu_1,\mu_2,\dots)$ such that
$\sum_i \mu_i = n$. The {\em length} $\ell(\mu)$ of $\mu$ is the
maximal index $i$ with $\mu_i$ nonzero. If, in addition, $\mu_1
\geq \mu_2 \geq \cdots$ then $\mu$ is also called a {\em
partition} and one writes $\mu \vdash n$ and $|\mu|=n$.

It is a standard result that the centralizer algebras of
$\slN$-modules are canonically subalgebras of $\C\Sigma_n$. So it
is not surprising to see that the centralizer algebras of
$\UslN$-modules are realized as subalgebras of $\CH_n(q)$, the
quantum deformation of $\C\Sigma_n$.

\medskip

Now fix $\g=\slN$ and let $V$ be the module of the fundamental
representation of $\UslN$. With suitable basis $\{ v_1,\dots,v_N
\}$ of $V$ and generators $\{ K_i^{\pm1}, E_i, F_i \mid 1 \leq i
\leq N-1 \}$ of $\UslN$, the fundamental representation is given
by the matrices
\begin{equation}
  \begin{array}{lll}
   K_i &\mapsto& q E_{ii} + q^{-1} E_{i+1,i+1} + \sum_{j \neq i} E_{jj}, \\
   E_i &\mapsto& E_{i,i+1}, \\
   F_i &\mapsto& E_{i+1,i},
  \end{array}
\end{equation}
where $E_{ij}$ is the $N \times N$ matrix with $1$ in the
$(i,j)$-position and $0$ elsewhere. We also have
\begin{eqnarray}
  && q^{1/N}\theta_V = q^N \cdot \id_V, \\
  \label{eqn:rho}
  && K_{2\rho}(v_i) = q^{N+1-2i} v_i,
\end{eqnarray}
and
\begin{equation}
  q^{1/N}\cR_{V,V}(v_i \otimes v_j) =
  \left\{ \begin{array}{ll}
    q v_i \otimes v_j, & i = j, \\
    v_j \otimes v_i, & i < j, \\
    v_j \otimes v_i + (q-q^{-1}) v_i \otimes v_j, & i > j. \\
  \end{array} \right.
\end{equation}

It is straightforward to verify that the homomorphism $h_V : \C
B_n \to \CC_n(V)$ factors through $\CH_n(q)$ via
\begin{equation}
  q^{1/N}\sigma_i \mapsto g_i \mapsto q^{1/N}h_V(\sigma_i).
\end{equation}
Therefore, $\Vn$ is a module of both $\UslN$ and $\CH_n(q)$, and
the two algebras act commutatively on $\Vn$. For convenience, we
introduce an $N$-independent homomorphism
\begin{equation}
  h : \C B_n \to \CH_n(q), \quad \sigma_i \mapsto g_i.
\end{equation}

\medskip

Let $S^\lambda$ denote the irreducible module of $\CH_n(q)$
indexed by the partition $\lambda \vdash n$ and let
$\zeta^\lambda$ denote its character. Fix a minimal projection
$p_\lambda \in \CH_\lambda(q)$ for each $\lambda \vdash n$. Let
$V_\lambda$ denote the irreducible $\UslN$-module, whose highest
weight vector $v$ behaves like $K_i(v) =
q^{\lambda_i-\lambda_{i+1}} v$, if $\ell(\lambda) \leq N$ and be
$0$ otherwise.

We state below two important facts concerning the $\UslN$-module
$\Vn$. One is the irreducible decomposition of $\UslN$-module
\begin{equation}\label{eqn:decom_svn}
  \Vn = \bigoplus_{\lambda \vdash n, \; \ell(\lambda) \leq N}
  \dim S^\lambda \cdot V_\lambda.
\end{equation}
in which the subspace $\dim S^\lambda \cdot V_\lambda$ is
$\CH_n(q)$-invariant and, as a $\CH_n(q)$-module, consists of only
$S^\lambda$-components. Notice that the $\UslN$-modules $\{
V_\lambda \mid \lambda \vdash n, \; \ell(\lambda) \leq N \}$ are
mutually inequivalent. Comparing (\ref{eqn:decom_svn}),
(\ref{eqn:decom_hn}) with (\ref{eqn:decom_vn}),
(\ref{eqn:decom_cn}), we have immediately
\begin{equation}
  \CC_n(V) = \bigoplus_{ \lambda \vdash n, \; \ell(\lambda) \leq N }
  \CH_\lambda(q).
\end{equation}
Moreover, for every partition $\lambda \vdash n$,
\begin{equation}
  p_\lambda \Vn \cong V_\lambda.
\end{equation}

The other fact is the weight decomposition of $\UslN$-module
\begin{equation}
  \Vn = \bigoplus_{\mu \models n, \; \ell(\mu) \leq N} M^\mu
\end{equation}
where
\begin{equation}
  M^\mu = \{ v \in \Vn \mid K_i(v) = q^{\mu_i-\mu_{i+1}}v \}.
\end{equation}
Moreover, the dimensions of the weight spaces of $V_\lambda$ for
$\lambda \vdash n$
\begin{equation}
  K_{\lambda\mu}
  = \dim (p_\lambda\Vn \cap M^\mu)
\end{equation}
are encoded in Schur polynomial as
\begin{equation}\label{eqn:Kostka}
  s_\lambda(z_1,\dots,z_N) = \sum_{\mu \models n, \; \ell(\mu) \leq N}
  K_{\lambda\mu} \cdot \prod_{j=1}^N z_j^{\mu_j}.
\end{equation}
Indeed, $M^\mu$ is nothing but the subspace of $\Vn$ spanned by
the vectors $v_{i_1} \otimes \cdots \otimes v_{i_n}$ in which
$v_i$ appears precisely $\mu_i$ times. It is clear that $M^\mu$ is
$\CH_n(q)$-invariant. In the literature, $M^\mu$ is called {\em
permutation module} and the integers $K_{\lambda\mu}$ are referred
to as {\em Kostka numbers}.

\medskip

Various choices of minimal projections of the Hecke algebras are
available in \cite{AM}\cite{DJ2}\cite{Gyo}\cite{Mur}. It is also
shown

\begin{thm}[Aiston-Morton {\cite[Theorem 5.5]{AM}}]
\label{thm:theta} For each partition $\lambda \vdash n$ with
$\ell(\lambda) \leq N$, one has
\begin{equation}
  \theta_{V_\lambda} = q^{\kappa_\lambda+nN-n^2/N} \cdot \id_{V_\lambda}
\end{equation}
where
\begin{equation}
  \kappa_\lambda
  = \sum_{i=1}^{\ell(\lambda)} \sum_{j=1}^{\lambda_i} 2(j-i).
\end{equation}
\end{thm}

The next proposition is a strong version of Lemma \ref{lem:trace}.
Equation (\ref{eqn:hecketrace}) holds even for $x \not\in
\CC_n(V)$.

\begin{prop}\label{prop:trace}
We have
\begin{equation}
  \dim_q V_\lambda = s_\lambda \Big( q^{N-1},q^{N-3},\dots,q^{-(N-1)} \Big)
\end{equation}
thus for every $x \in \CH_n(q)$,
\begin{equation}\label{eqn:hecketrace}
  \tr_\Vn x = \sum_{\lambda \vdash n} \zeta^\lambda(x) \cdot
  s_\lambda \Big( q^{N-1},q^{N-3},\dots,q^{-(N-1)} \Big).
\end{equation}
\end{prop}

\begin{proof}
By (\ref{eqn:k2rho}) and (\ref{eqn:rho}), $K_{2\rho}$ acts as a
scalar $\prod_{i=1}^{N} q^{(N+1-2i)\mu_j}$ on $M^\mu$. Therefore,
it follows from identity (\ref{eqn:Kostka}) that for each $\lambda
\vdash n$,
\begin{equation}
\begin{split}
  \dim_q V_\lambda
  = \sum_{\mu \models n, \; \ell(\mu) \leq N}
  K_{\lambda\mu} \cdot \prod_{i=1}^{N} q^{(N+1-2i)\mu_j}
  = s_\lambda \Big( q^{N-1},q^{N-3},\dots,q^{-(N-1)} \Big).
\end{split}
\end{equation}
\end{proof}

Now it is time to give our main result.

\begin{thm}\label{thm:inv}
Let $\CL$ be an oriented link with $l$ components
$\CL_1,\dots,\CL_l$. Suppose $\CL$ is the closure of $\beta \in
B_m$ and the $m$ strands of $\beta$ are living on
$\CL_{i_1},\dots,\CL_{i_m}$, respectively. Then for partitions
$\lambda^i \vdash n_i$, $i=1,\dots,l$, we have
\begin{equation}\label{eqn:inv}
\begin{split}
  I_{\CL;V_{\lambda^1},\dots,V_{\lambda^l}}
  &= q^{-\sum_{i=1}^l (\kappa_{\lambda^i}+n_iN-n_i^2/N) w(\CL_i) -
  w(\beta^{(n_{i_1},\dots,n_{i_m})})/N} \cdot \\
  & \sum_{\lambda \vdash n}
  \zeta^\lambda(x) \cdot
  s_\lambda \Big( q^{N-1},q^{N-3},\dots,q^{-(N-1)} \Big),
\end{split}
\end{equation}
where $n=n_{i_1}+\cdots+n_{i_m}$, $\beta^{(n_{i_1},\dots,n_{i_m})}
\in B_n$ is the braid obtained by cabling the $j$-th strand of
$\beta$ to $n_{i_j}$ parallel ones and $x =
h(\beta^{(n_{i_1},\dots,n_{i_m})}) \cdot (p_{\lambda^{i_1}}
\otimes \cdots \otimes p_{\lambda^{i_m}}) \in \CH_n(q)$.
\end{thm}

\begin{proof}
Combine Lemma \ref{lem:cable}, Theorem \ref{thm:theta} and
Proposition \ref{prop:trace}.
\end{proof}

One notices that, on the right hand side of (\ref{eqn:inv}), there
is an explicit factor $q^{1/N}$ to the power
\begin{equation}
  \sum_{i=1}^l n_i^2 w(\CL_i) - w(\beta^{(n_{i_1},\dots,n_{i_m})})
  = -2\sum_{i<j} n_in_j \lk(\CL_i,\CL_j)
\end{equation}
where $\lk(\CL_i,\CL_j)$ are the linking numbers. As in
\cite{Labastida-Marino2}, we drop this insignificant factor and
regard the remaining part as a rational function of $q$ and $q^N$.

\begin{defn}
The {\em colored HOMFLY polynomial}
$W_{\CL;\lambda^1,\dots,\lambda^l}(t,\nu)$ with $\lambda^i \vdash
n_i$ is a rational function of $t^{1/2}, \nu^{1/2}$ determined by
\begin{equation}
  W_{\CL;\lambda^1,\dots,\lambda^l}(t,\nu)
  |_{t^{1/2}=q^{-1}, \; \nu^{1/2}=q^{-N}}
  = q^{2\sum_{i<j} n_in_j \lk(\CL_i,\CL_j)/N} \cdot
  I_{\CL;V_{\lambda^1},\dots,V_{\lambda^l}}.
\end{equation}
\end{defn}

Note that the definition means the components of the link $\CL$
are labeled by partitions rather than $\UslN$-modules. When the
labeling partitions are trivial (the unique partition of $1$), the
colored HOMFLY polynomial, up to a simple factor, specializes to
the HOMFLY polynomial:
\begin{equation}
  P_\CL(t,\nu)
  = \nu^{\lk(\CL)} \cdot
  \frac {t^{1/2}-t^{-1/2}} {\nu^{1/2}-\nu^{-1/2}} \cdot
  W_{\CL;(1),\dots,(1)}(t,\nu)
\end{equation}

\medskip

Let $s^*_\lambda(t,\nu)$ be defined by (see (\ref{eqn:f}))
\begin{equation}
  s^*_\lambda(t,\nu)
  |_{t^{1/2}=q^{-1}, \; \nu^{1/2}=q^{-N}}
  = s_\lambda \Big( q^{N-1},q^{N-3},\dots,q^{-(N-1)} \Big).
\end{equation}

\begin{cor}\label{cor:inv}
In the same notations as Theorem \ref{thm:inv}, we have
\begin{equation}
  W_{\CL;\lambda^1,\dots,\lambda^l} (t,\nu)
  = t^{\sum_{i=1}^l \kappa_{\lambda^i} w(\CL_i)/2} \cdot
  \nu^{\sum_{i=1}^l n_i w(\CL_i)/2} \cdot
  \sum_{\lambda \vdash n} \zeta^\lambda(x) \Big|_{q=t^{-1/2}} \cdot
  s^*_\lambda(t,\nu).
\end{equation}
\end{cor}

\section{Torus links}\label{sec:torus}

Let notations be the same as in the previous section. In this
section, we derive an explicit formula of the colored HOMFLY
polynomial of torus links by applying Corollary \ref{cor:inv}.

The {\em torus link} $T(r,k)$ is defined to be the closure of
$(\delta_r)^k = (\sigma_1\cdots\sigma_{r-1})^k$. They form the
family of link that can be put on the standardly embedded torus $T
\subset \R^3$. Some common links such as the trefoil knot
$T(2,3)$, the Hopf link $T(2,2)$ are included in this family.

\begin{thm}\label{thm:torus}
Let $\CL$ be the torus link $T(rl,kl)$ with $r,k$ relatively
prime. Let $\lambda^i \vdash n_i$, $i=1,\dots,l$ be partitions and
$n=n_1+\cdots+n_l$. Then
\begin{equation}
  W_{\CL;\lambda^1,\dots,\lambda^l}(t,\nu)
  = t^{kr\sum_{i=1}^l \kappa_{\lambda^i}/2} \cdot
  \nu^{k(r-1)n/2} \cdot
  \sum_{\lambda \vdash rn}
  c^\lambda_{\lambda^1\dots\lambda^l} \cdot
  t^{-k\kappa_\lambda/2r} \cdot s^*_\lambda(t,\nu)
\end{equation}
where $c^\lambda_{\lambda^1\dots\lambda^l}$ are the integers
determined by the equation
\begin{equation}
  \prod_{i=1}^l s_{\lambda^i}(x_1^r,x_2^r,\dots)
  = \sum_{\lambda \vdash rn}
  c^\lambda_{\lambda^1\dots\lambda^l} \cdot
  s_\lambda(x_1,x_2,\dots).
\end{equation}
\end{thm}

The theorem is an easy consequence of following lemmas.

\begin{lem}\label{lem:twist}
For each partition $\lambda \vdash n$ we have
\begin{equation}
  h(\Delta_n^2) \cdot p_\lambda
  = q^{\kappa_\lambda} \cdot p_\lambda.
\end{equation}
\end{lem}

\begin{proof}
Compare Lemma \ref{lem:theta} with Theorem \ref{thm:theta}.
\end{proof}

\begin{lem}\label{lem:delta}
Let $\lambda^i \vdash n_i$, $i=1,\dots,l$ be partitions and
$n=n_1+\cdots+n_l$. Let $r, k$ be relatively prime integers and
$\beta \in B_{rn}$ be the braid obtained by cabling the
$(il+j)$-th strand of $(\delta_{rl})^{kl}$ to $n_j$ parallel ones.
Then, for each partition $\lambda \vdash rn$ we have
\begin{equation}\label{eqn:delta}
  \zeta^\lambda \Big( h(\beta) \cdot
  (p_{\lambda^1} \otimes \cdots \otimes p_{\lambda^l})^\otr \Big)
  = c^\lambda_{\lambda^1\dots\lambda^l} \cdot
  q^{-k \sum_{i=1}^l \kappa_{\lambda^i} + k\kappa_\lambda/r}.
\end{equation}
\end{lem}

\begin{proof}
Put $p = p_{\lambda^1} \otimes \cdots \otimes p_{\lambda^l}$ and
let $\pi_\lambda$ be the unit of $\CH_\lambda(q)$. Note that
$\pi_\lambda$ is a central element of $\CH_{rn}(q)$ and $h(\beta)$
is commutative with $p^\otr$. Then
\begin{equation}
  x_\lambda = \pi_\lambda \cdot h(\beta) \cdot p^\otr
\end{equation}
is a matrix in $\CH_\lambda(q)$, whose normal trace is
\begin{equation}
  \tr x_\lambda
  = \zeta^\lambda \Big( h(\beta) \cdot p^\otr \Big).
\end{equation}
Notice that
\begin{equation}
  h(\beta^r)
  = h(\Delta_{rn}^{2k}) \cdot
  \Big( h(\Delta_{n_1}^{-2k}) \otimes \cdots \otimes
  h(\Delta_{n_l}^{-2k}) \Big)^\otr.
\end{equation}
It follows from Lemma \ref{lem:twist} that
\begin{equation}
  x_\lambda^r = \pi_\lambda \cdot h(\beta^r) \cdot p^\otr
  = q^{-kr \sum_{i=1}^l \kappa_{\lambda^i} + k\kappa_\lambda} \cdot
  \pi_\lambda \cdot p^\otr.
\end{equation}
Therefore, the eigenvalues of $x_\lambda$ are either $0$ or $q^{-k
\sum_{i=1}^l \kappa_{\lambda^i} + k\kappa_\lambda/r}$ times an
$r$-th root of unity, for $\pi_\lambda \cdot p^\otr \in
\CH_\lambda(q)$ is also a projection. Since $\tr x_\lambda$ is
always a rational function of $q$ in rational coefficients (easily
seen with suitable choice of minimal projections), it follows that
$\tr x_\lambda$ has to be $q^{-k \sum_{i=1}^l \kappa_{\lambda^i} +
k\kappa_\lambda/r}$ times a rational number $a^\lambda$ which is
independent of $q$.

Now let $q \to 1$. Passing to the limit, $h(\beta)$ degenerates to
a permutation $\tau \in \Sigma_{rn}$ which acts cyclicly on the
$\Vn$-factors of $\V{rn} = \Vn \otimes \cdots \otimes \Vn$,
because $r,k$ are relatively prime. By identity
(\ref{eqn:Kostka}), we have
\begin{equation*}
\begin{split}
  & \sum_{\lambda \vdash rn} a^\lambda \cdot
  s_\lambda(z_1,\dots,z_N) \\
  =& \sum_{\lambda \vdash rn} a^\lambda
  \sum_{\mu \models rn, \; \ell(\mu) \leq N}
  \dim (p_\lambda\V{rn} \cap M^\mu) \cdot
  \prod_{j=1}^N z_j^{\mu_j} \\
  =& \sum_{\mu \models rn, \; \ell(\mu) \leq N}
  \tr \tau|_{p^\otr\V{rn} \cap M^\mu} \cdot
  \prod_{j=1}^N z_j^{\mu_j} \\
  =& \sum_{\mu \models n, \; \ell(\mu) \leq N}
  \dim (p\Vn \cap M^\mu) \cdot
  \prod_{j=1}^N z_j^{r\mu_j} \\
  =& \prod_{i=1}^l
  \Big( \sum_{\mu \models n_i, \; \ell(\mu) \leq N}
  \dim (p_{\lambda^i}\V{n_i} \cap M^\mu) \cdot
  \prod_{j=1}^N z_j^{r\mu_j} \Big) \\
  =& \prod_{i=1}^l s_{\lambda^i}(z_1^r,\dots,z_N^r).
\end{split}
\end{equation*}
Since above equality holds for all $N$, we must have $a^\lambda =
c^\lambda_{\lambda^1\dots\lambda^l}$.
\end{proof}

\begin{rem}
In the case $l=2$ and $r=1$, $k=0$, equation (\ref{eqn:delta})
specializes to
\begin{equation}
  \zeta^\lambda (p_{\lambda^1} \otimes p_{\lambda^2})
  = c^\lambda_{\lambda^1\lambda^2}
\end{equation}
which implies that (see (\ref{eqn:proj_decom}))
\begin{equation}
  V_{\lambda^1} \otimes V_{\lambda^2} =
  \bigoplus_{\lambda \vdash n}
  c^\lambda_{\lambda^1\lambda^2} \cdot V_\lambda.
\end{equation}
In the literature, the integers $c^\lambda_{\lambda^1\lambda^2}$
are referred to as {\em Littlewood-Richardson coefficients}.
\end{rem}

\begin{proof}[Proof of Theorem 5.1]
By Corollary \ref{cor:inv} and Lemma \ref{lem:delta}, we have
\begin{equation*}
\begin{split}
  & W_{\CL;\lambda^1,\dots,\lambda^l}(t,\nu) \\
  = & t^{\sum_{i=1}^l \kappa_{\lambda^i} k(r-1)/2} \cdot
  \nu^{\sum_{i=1}^l n_i k(r-1)/2} \cdot
  \sum_{\lambda \vdash rn}
  c^\lambda_{\lambda^1\dots\lambda^l} \cdot
  t^{k\sum_{i=1}^l \kappa_{\lambda^i}/2 - k\kappa_\lambda/2r} \cdot
  s^*_\lambda(t,\nu) \\
  = & t^{kr\sum_{i=1}^l \kappa_{\lambda^i}/2} \cdot
  \nu^{k(r-1)n/2} \cdot
  \sum_{\lambda \vdash rn}
  c^\lambda_{\lambda^1\dots\lambda^l} \cdot
  t^{-k\kappa_\lambda/2r} \cdot
  s^*_\lambda(t,\nu).
\end{split}
\end{equation*}
\end{proof}

The functions $s^*_\lambda(t,\nu)$ and the coefficients
$c^\lambda_{\lambda^1\dots\lambda^l}$ can be computed by using the
Frobenius formula as follows. Let $\chi^\lambda$ and $C_\mu$
denote the character and conjugacy class of the symmetric group
$\Sigma_n$ indexed by $\lambda, \mu \vdash n$. The Frobenius
formula says the Newton polynomial
\begin{equation}
  p^\mu(x_1,x_2,\dots)
  = \prod_{i=1}^{\ell(\mu)} \sum_{j\geq1} x_j^{\mu_i}
\end{equation}
is expressed in terms of Schur polynomials as
\begin{equation}
  p^\mu(x_1,x_2,\dots)
  = \sum_{\lambda \vdash |\mu|}
  \chi^\lambda(C_\mu) \cdot
  s_\lambda(x_1,x_2,\dots).
\end{equation}
Its inverse for $\lambda \vdash n$ is
\begin{equation}
  s_\lambda(x_1,x_2,\dots)
  = \sum_{\mu \vdash n}
  \frac{|C_\mu|}{n!} \chi^\lambda(C_\mu) \cdot
  p^\mu(x_1,x_2,\dots).
\end{equation}
Therefore, for partition $\lambda \vdash n$ we have
\begin{equation}\label{eqn:f}
  s^*_\lambda(t,\nu)
  = \sum_{\mu \vdash n}
  \frac{|C_\mu|}{n!} \chi^\lambda(C_\mu)
  \prod_{i=1}^{\ell(\mu)}
  \frac {\nu^{\mu_i/2}-\nu^{-\mu_i/2}} {t^{\mu_i/2}-t^{-\mu_i/2}}.
\end{equation}
Moreover, it is clear that
\begin{equation}
  p^{\mu^1}(x_1,x_2,\dots) \cdot p^{\mu^2}(x_1,x_2,\dots)
  = p^{\mu^1+\mu^2}(x_1,x_2,\dots)
\end{equation}
and
\begin{equation}
  p^\mu(x_1^r,x_2^r,\dots)
  = p^{\mu_{(r)}}(x_1,x_2,\dots)
\end{equation}
where $\mu^1+\mu^2$ is the partition in which the number of each
positive integer is the sum of those in $\mu^1, \mu^2$ and
$\mu_{(r)}$ means the partition $(r\mu_1,r\mu_2,\dots)$. Hence,
for partitions $\lambda^i \vdash n_i$, $i=1,\dots,l$ and $\lambda
\vdash r(n_1+\cdots+n_l)$,
\begin{equation}
  c^\lambda_{\lambda^1\dots\lambda^l}
  =\sum_{\mu^1 \vdash n_1}
  \frac{|C_{\mu^1}|}{n_1!} \chi^{\lambda^1}(C_{\mu^1}) \cdots
  \sum_{\mu^l \vdash n_l}
  \frac{|C_{\mu^l}|}{n_l!} \chi^{\lambda^l}(C_{\mu^l}) \cdot
  \chi^\lambda( C_{(\mu^1+\cdots+\mu^l)_{(r)}} ).
\end{equation}

We finish this section by offering the following sample
calculations.

\begin{exam}
Torus knot $T(2,k)$, $k \not\equiv 0 \pmod 2$.
\begin{equation}
\begin{split}
  & W_{(1)}(t,\nu)
  = \nu^{k/2} \Big(
    t^{-k/2} s^*_{(2)}(t,\nu)
    - t^{k/2} s^*_{(1,1)}(t,\nu)
  \Big), \\
  & W_{(2)}(t,\nu)
  = \nu^{k} \Big(
    t^{-k} s^*_{(4)}(t,\nu)
    - t^{k} s^*_{(3,1)}(t,\nu)
    + t^{2k} s^*_{(2,2)}(t,\nu)
  \Big), \\
  & W_{(1,1)}(t,\nu)
  = \nu^{k} \Big(
    t^{-2k} s^*_{(2,2)}(t,\nu)
    - t^{-k} s^*_{(2,1,1)}(t,\nu)
    + t^{k} s^*_{(1,1,1,1)}(t,\nu)
  \Big), \\
  & W_{(3)}(t,\nu)
  = \nu^{3k/2} \Big(
    t^{-3k/2} s^*_{(6)}(t,\nu)
    - t^{3k/2} s^*_{(5,1)}(t,\nu)
  \\ & \quad \quad \quad \quad \quad
    + t^{7k/2} s^*_{(4,2)}(t,\nu)
    - t^{9k/2} s^*_{(3,3)}(t,\nu)
  \Big), \\
  & W_{(2,1)}(t,\nu)
  = \nu^{3k/2} \Big(
    t^{-5k/2} s^*_{(4,2)}(t,\nu)
    - t^{-3k/2} s^*_{(4,1,1)}(t,\nu)
    - t^{-3k/2} s^*_{(3,3)}(t,\nu)
  \\ & \quad \quad \quad \quad \quad
    + t^{3k/2} s^*_{(2,2,2)}(t,\nu)
    + t^{3k/2} s^*_{(3,1,1,1)}(t,\nu)
    - t^{5k/2} s^*_{(2,2,1,1)}(t,\nu)
  \Big), \\
  & W_{(1,1,1)}(t,\nu)
  = \nu^{3k/2} \Big(
    t^{-9k/2} s^*_{(2,2,2)}(t,\nu)
    - t^{-7k/2} s^*_{(2,2,1,1)}(t,\nu)
  \\ & \quad \quad \quad \quad \quad
    + t^{-3k/2} s^*_{(2,1,1,1,1)}(t,\nu)
    - t^{3k/2} s^*_{(1,1,1,1,1,1)}(t,\nu)
  \Big).
\end{split}
\end{equation}
In particular,
\begin{equation}
  W_{(1)}(t,\nu)
  = \frac {\nu^{1/2}-\nu^{-1/2}} {t^{1/2}-t^{-1/2}}
  \Big(
    \frac {t^{\frac{k+1}{2}}-t^{-\frac{k+1}{2}}} {t-t^{-1}} \nu^{\frac{k-1}{2}}
    - \frac {t^{\frac{k-1}{2}}-t^{-\frac{k-1}{2}}} {t-t^{-1}} \nu^{\frac{k+1}{2}}
  \Big).
\end{equation}
\end{exam}

\begin{exam}
Torus knot $T(3,k)$, $k \not\equiv 0 \pmod 3$.
\begin{equation}
\begin{split}
  & W_{(1)}(t,\nu)
  = \nu^{k} \Big(
    t^{-k} s^*_{(3)}(t,\nu)
    - s^*_{(2,1)}(t,\nu)
    + t^{k} s^*_{(1,1,1)}(t,\nu)
  \Big), \\
  & W_{(2)}(t,\nu)
  = \nu^{2k} \Big(
    t^{-2k} s^*_{(6)}(t,\nu)
    - s^*_{(5,1)}(t,\nu)
    + t^{2k} s^*_{(4,1,1)}(t,\nu)
  \\ & \quad \quad \quad \quad \quad
    + t^{2k} s^*_{(3,3)}(t,\nu)
    - t^{3k} s^*_{(3,2,1)}(t,\nu)
    + t^{4k} s^*_{(2,2,2)}(t,\nu)
  \Big), \\
  & W_{(1,1)}(t,\nu)
  = \nu^{2k} \Big(
    t^{-4k} s^*_{(3,3)}(t,\nu)
    - t^{-3k} s^*_{(3,2,1)}(t,\nu)
    + t^{-2k} s^*_{(3,1,1,1)}(t,\nu)
  \\ & \quad \quad \quad \quad \quad
    + t^{-2k} s^*_{(2,2,2)}(t,\nu)
    - s^*_{(2,1,1,1,1)}(t,\nu)
    + t^{2k} s^*_{(1,1,1,1,1,1)}(t,\nu)
  \Big). \\
\end{split}
\end{equation}
\end{exam}

\begin{exam}
Torus link $T(2,2k)$.
\begin{equation}
\begin{split}
  & W_{(1),(1)}(t,\nu)
  = t^{-k} s^*_{(2)}(t,\nu)
    + t^{k} s^*_{(1,1)}(t,\nu), \\
  & W_{(2),(1)}(t,\nu)
  = t^{-2k} s^*_{(3)}(t,\nu)
    + t^{k} s^*_{(2,1)}(t,\nu), \\
  & W_{(1,1),(1)}(t,\nu)
  = t^{-k} s^*_{(2,1)}(t,\nu)
    + t^{2k} s^*_{(1,1,1)}(t,\nu), \\
  & W_{(2),(2)}(t,\nu)
  = t^{-4k} s^*_{(4)}(t,\nu)
    + s^*_{(3,1)}(t,\nu)
    + t^{2k} s^*_{(2,2)}(t,\nu), \\
  & W_{(2),(1,1)}(t,\nu)
  = t^{-2k} s^*_{(3,1)}(t,\nu)
    + t^{2k} s^*_{(2,1,1)}(t,\nu), \\
  & W_{(1,1),(1,1)}(t,\nu)
  = t^{-2k} s^*_{(2,2)}(t,\nu)
    + s^*_{(2,1,1)}(t,\nu)
    + t^{4k} s^*_{(1,1,1,1)}(t,\nu).
\end{split}
\end{equation}
In particular,
\begin{equation}
  W_{(1),(1)}(t,\nu)
  = \frac {\nu^{1/2}-\nu^{-1/2}} {t^{1/2}-t^{-1/2}}
  \Big(
    \frac {t^{\frac{2k-1}{2}}+t^{-\frac{2k-1}{2}}} {t-t^{-1}} \nu^{\frac12}
    - \frac {t^{\frac{2k+1}{2}}+t^{-\frac{2k+1}{2}}} {t-t^{-1}} \nu^{-\frac12}
  \Big).
\end{equation}
\end{exam}

\begin{exam}
Torus link $T(3,3k)$.
\begin{equation}
\begin{split}
  & W_{(1),(1),(1)}(t,\nu)
  = t^{-3k} s^*_{(3)}(t,\nu)
    + 2 s^*_{(2,1)}(t,\nu)
    + t^{3k} s^*_{(1,1,1)}(t,\nu), \\
  & W_{(2),(1),(1)}(t,\nu)
  = t^{-5k} s^*_{(4)}(t,\nu)
    + 2 t^{-k} s^*_{(3,1)}(t,\nu)
    + t^{k} s^*_{(2,2)}(t,\nu)
    + t^{3k} s^*_{(2,1,1)}(t,\nu), \\
  & W_{(1,1),(1),(1)}(t,\nu)
  = t^{-3k} s^*_{(3,1)}(t,\nu)
    + t^{-k} s^*_{(2,2)}(t,\nu)
    + 2 t^{k} s^*_{(2,1,1)}(t,\nu)
    + t^{5k} s^*_{(1,1,1,1)}(t,\nu).
\end{split}
\end{equation}
\end{exam}

\section{On the Labastida-Mari\~no-Vafa conjecture}

As before, we have an oriented link $\CL$ with $l$ components.
Define the generating function
\begin{equation}
  Z(\x_1,\dots,\x_l)
  = \sum_{\lambda^1,\dots,\lambda^l}
  W_{\lambda^1,\dots,\lambda^l}(t,\nu) \cdot
  s_{\lambda^1}(\x_1) \cdots s_{\lambda^l}(\x_l)
\end{equation}
where each $\x_i = \{x_{i,1},x_{i,2},\dots\}$ is a set of
indeterminate variables and $\lambda^i$ runs over all partitions
including the empty one (the unique partition of zero). When all
$\lambda^i$ are empty, the summand gives rise to the leading term
$1$.

One can expand $\log Z(\x_1,\dots,\x_l)$ as
\begin{equation}
  \log Z(\x_1,\dots,\x_l)
  = \sum_{d=1}^\infty \sum_{\lambda^1,\dots,\lambda^l}
  \frac{1}{d}\, f_{\lambda^1,\dots,\lambda^l}(t^d,\nu^d) \cdot
  s_{\lambda^1}(\x_1^d) \cdots s_{\lambda^l}(\x_l^d)
\end{equation}
where $\x_i^d = \{x_{i,1}^d,x_{i,2}^d,\dots\}$. See
\cite{Labastida-Marino2} for an explanation why such an expansion
exists by using the so-called plethystic exponential. The
functions $f_{\lambda^1,\dots,\lambda^l}(t,\nu)$ are referred to
as the {\em reformulated colored HOMFLY polynomial}.
Labastida-Mari\~no-Vafa conjecture says that these functions have
the following highly nontrivial structure.

Write for $\lambda,\mu \vdash n$
\begin{equation}
  M_{\lambda\mu}(t)
  = \sum_{\tau \vdash n} \frac{|C_\tau|}{n!}
  \chi^\lambda(C_\tau) \chi^\mu(C_\tau) \cdot
  \frac {\prod_{j=1}^{\ell(\tau)} (t^{-\tau_j/2}-t^{\tau_j/2})}
  {t^{-1/2}-t^{1/2}}.
\end{equation}

\begin{conj}[Labastida-Mari\~no-Vafa \cite{Labastida-Marino2}\cite{LMV}]
For partitions $\lambda^1,\dots,\lambda^l$,
\begin{eqnarray}
  \label{eqn:LM}
  && f_{\lambda^1,\dots,\lambda^l}(t,\nu)
  = \sum_{\mu^1 \vdash |\lambda^1|,\dots,\mu^l \vdash |\lambda^l|}
  \hat{f}_{\mu^1,\dots,\mu^l}(t,\nu) \cdot
  M_{\lambda^1\mu^1}(t) \cdots M_{\lambda^l\mu^l}(t), \\
  && \hat{f}_{\mu^1,\dots,\mu^l}(t,\nu)
  = \sum_{g \geq 0} \sum_Q
  N_{\mu^1,\dots,\mu^l,g,Q} \cdot (t^{1/2}-t^{-1/2})^{2g+l-2} \cdot \nu^Q,
\end{eqnarray}
where $N_{\mu^1,\dots,\mu^l,g,Q}$ are integers and $Q$ are either
all integers or all semi-integers.
\end{conj}

Moreover, the integers $N_{\mu^1,\dots,\mu^l,g,Q}$ are interpreted
as quantities involved in the enumerative geometry of the resolved
conifold. See \cite{Gopakumar-Vafa}\cite{Ooguri-Vafa} for example.

\medskip

Till now, besides the trivial links, the conjecture was verified
only for some simplest knots and links with small partitions. A
proof of the Labastida-Mari\~no-Vafa conjecture seems to appeal to
deep knowledge of mathematics and string theory. Using the formula
in the previous section, we can verify this conjecture for several
infinite families of torus links with small partitions. 
Our calculation also suggests 
a new structure of the reformulated colored
HOMFLY polynomial of torus links. Let us make this more precise first.

Define symmetric functions for $\lambda,\mu \vdash n$
\begin{equation}
  S_{\lambda,\mu}(\x)
  = \sum_{\tau \vdash n} \frac{|C_\tau|}{n!}
  \chi^\lambda(C_\tau) \chi^\mu(C_\tau) \cdot p^\tau(\x)
\end{equation}
and
\begin{equation}
\begin{split}
  s_{\mu;q}(\x)
  = & \sum_{\lambda \vdash n}
  (q-q^{-1}) M_{\lambda\mu}(q^{-2}) \cdot
  s_\lambda(\x) \\
  = & \sum_{\tau \vdash n} \frac{|C_\tau|}{n!}
  \chi^\mu(C_\tau) \cdot
  \prod_{j=1}^{\ell(\tau)} (q^{\tau_j}-q^{-\tau_j}) \cdot
  p^\tau(\x).
\end{split}
\end{equation}
We have the following conjecture for torus links.

\begin{conj}
For torus link $T(rl,kl)$ with $r,k$ relatively prime and
$n=n_1+\cdots+n_l$,
\begin{equation}\label{eqn:LMX}
\begin{split}
  & \sum_{\lambda^1 \vdash n_1,\dots,\lambda^l \vdash n_l}
  f_{\lambda^1,\dots,\lambda^l}(t,\nu) \cdot
  s_{\lambda^1}(\x_1) \cdots s_{\lambda^l}(\x_l)
  |_{t^{1/2}=q^{-1}, \; \nu^{1/2}=q^{-N}} \\
  = & \sum_{\lambda^1 \vdash n_1,\dots,\lambda^l \vdash n_l}
  (q^{k}-q^{-k})^{-2} \cdot q^{-k(r-1)nN} \cdot
  \sum_{\lambda \vdash rn}
  g_{\lambda^1,\dots,\lambda^l}^\lambda(q^{2k}) \cdot \\
  & s_{\lambda;q^{k}} \Big( q^{N-1},q^{N-3},\dots,q^{-(N-1)} \Big) \cdot
  s_{\lambda^1;q^{k}}(\x_1) \cdots s_{\lambda^l;q^{k}}(\x_l).
\end{split}
\end{equation}
where $g_{\lambda^1,\dots,\lambda^l}^\lambda(t) \in \Z[t^{\pm1}]$
is invariant under $t \to t^{-1}$.
\end{conj}

The following examples verify Conjecture 6.2.

\begin{exam}
Torus knot $T(2,k)$, $k \not\equiv 0 \pmod 2$. Followings are
nonvanishing $g_{\lambda'}^\lambda(t)$'s for $|\lambda'| \leq 4$.
\begin{equation}
\begin{split}
  & g_{(1)}^{(2)}(t) = 1, \\
  & g_{(1,1)}^{(2,2)}(t) = 1, \\
  & g_{(1,1,1)}^{(2,2,2)}(t) = t+t^{-1}, \\
  & g_{(2,1)}^{(2,2,2)}(t) =
  g_{(1,1,1)}^{(3,2,1)}(t) = 1, \\
  & g_{(1,1,1,1)}^{(2,2,2,2)}(t) = t^3+2t+1+2t^{-1}+t^{-3}, \\
  & g_{(2,1,1)}^{(2,2,2,2)}(t) =
  g_{(1,1,1,1)}^{(3,2,2,1)}(t) = t^2+t+2+t^{-1}+t^{-2}, \\
  & g_{(2,1,1)}^{(3,2,2,1)}(t) = t+1+t^{-1}, \\
  & g_{(2,2)}^{(2,2,2,2)}(t) =
  g_{(1,1,1,1)}^{(4,2,2)}(t) =
  g_{(1,1,1,1)}^{(3,3,1,1)}(t) = t+t^{-1}, \\
  & g_{(3,1)}^{(2,2,2,2)}(t) =
  g_{(2,2)}^{(3,2,2,1)}(t) =
  g_{(2,1,1)}^{(4,2,2)}(t) =
  g_{(2,1,1)}^{(3,3,1,1)}(t) = 1, \\
  & g_{(1,1,1,1)}^{(4,3,1)}(t) =
  g_{(1,1,1,1)}^{(4,2,1,1)}(t) =
  g_{(1,1,1,1)}^{(3,3,2)}(t) = 1.
\end{split}
\end{equation}
\end{exam}

\begin{exam}
Torus knot $T(3,k)$, $k \not\equiv 0 \pmod 3$. Followings are
nonvanishing $g_{\lambda'}^\lambda(t)$'s for $|\lambda'| \leq 3$.
\begin{equation}
\begin{split}
  & g_{(1)}^{(3)}(t) = 1, \\
  & g_{(1,1)}^{(3,3)}(t) = t+t^{-1}, \\
  & g_{(2)}^{(3,3)}(t) =
  g_{(1,1)}^{(4,2)}(t) = 1, \\
  & g_{(1,1,1)}^{(3,3,3)}(t) = t^4+2t^2+2t+2+2t^{-1}+2t^{-2}+t^{-4}, \\
  & g_{(2,1)}^{(3,3,3)}(t) =
  g_{(1,1,1)}^{(4,3,2)}(t) = t^3+t^2+2t+3+2t^{-1}+t^{-2}+t^{-3}, \\
  & g_{(2,1)}^{(4,3,2)}(t) = t^2+2t+2+2t^{-1}+t^{-2}, \\
  & g_{(1,1,1)}^{(5,3,1)}(t) = t^2+t+2+t^{-1}+t^{-2}, \\
  & g_{(2,1)}^{(5,3,1)}(t) = t+1+t^{-1}, \\
  & g_{(3)}^{(3,3,3)}(t) =
  g_{(1,1,1)}^{(6,3)}(t) =
  g_{(1,1,1)}^{(5,2,2)}(t) =
  g_{(1,1,1)}^{(4,4,1)}(t) = t+t^{-1}, \\
  & g_{(3)}^{(4,3,2)}(t) =
  g_{(2,1)}^{(6,3)}(t) =
  g_{(2,1)}^{(5,2,2)}(t) =
  g_{(2,1)}^{(4,4,1)}(t) =
  g_{(1,1,1)}^{(6,2,1)}(t) =
  g_{(1,1,1)}^{(5,4)}(t) = 1.
\end{split}
\end{equation}
\end{exam}

\begin{exam}
Torus link $T(2,2k)$. Followings are nonvanishing
$g_{\lambda^1,\lambda^2}^\lambda(t)$'s for $|\lambda^1| +
|\lambda^2| \leq 5$, up to symmetry
$g_{\lambda^1,\lambda^2}^\lambda(t) =
g_{\lambda^2,\lambda^1}^\lambda(t)$.
\begin{equation}
\begin{split}
  & g_{(0),(1)}^{(1)}(t) =
  g_{(1),(1)}^{(2)}(t) =
  g_{(2),(1)}^{(3)}(t) =
  g_{(3),(1)}^{(4)}(t) =
  g_{(4),(1)}^{(5)}(t) = 1, \\
  & g_{(2),(2)}^{(4)}(t) = t+1+t^{-1}, \\
  & g_{(2),(1,1)}^{(4)}(t) =
  g_{(2),(2)}^{(3,1)}(t) = 1, \\
  & g_{(3),(2)}^{(5)}(t) = t^2+t+3+t^{-1}+t^{-2}, \\
  & g_{(3),(2)}^{(4,1)}(t) =
  g_{(3),(1,1)}^{(5)}(t) =
  g_{(2,1),(2)}^{(5)}(t) = t+1+t^{-1}, \\
  & g_{(3),(2)}^{(3,2)}(t) =
  g_{(3),(1,1)}^{(4,1)}(t) =
  g_{(2,1),(2)}^{(4,1)}(t) =
  g_{(2,1),(1,1)}^{(5)}(t) = 1.
\end{split}
\end{equation}
\end{exam}

\medskip

One notices that (\ref{eqn:LM}) is indeed equivalent to
\begin{equation}
\begin{split}
  & \sum_{\lambda^1 \vdash n_1,\dots,\lambda^l \vdash n_l}
  f_{\lambda^1,\dots,\lambda^l}(t,\nu) \cdot
  s_{\lambda^1}(\x_1) \cdots s_{\lambda^l}(\x_l) \\
  = & \sum_{\mu^1 \vdash n_1,\dots,\mu^l \vdash n_l}
  (t^{-1/2}-t^{1/2})^{-l} \cdot
  \hat{f}_{\mu^1,\dots,\mu^l}(t,\nu) \cdot
  s_{\mu^1;t^{-1/2}}(\x_1) \cdots s_{\mu^l;t^{-1/2}}(\x_l).
\end{split}
\end{equation}
Therefore, it follows from the identities
\begin{eqnarray}
  && s_{\lambda;q^k} \Big( q^{N-1},q^{N-3},\dots,q^{-(N-1)} \Big)
  = s_{\lambda;q^N} \Big( q^{k-1},q^{k-3},\dots,q^{-(k-1)} \Big), \\
  && s_{\lambda;q^k}(\x)
  = \sum_{\mu \vdash |\lambda|}
  S_{\lambda,\mu} \Big( q^{k-1},q^{k-3},\dots,q^{-(k-1)} \Big) \cdot
  s_{\mu;q}(\x)
\end{eqnarray}
that (\ref{eqn:LMX}) implies
\begin{equation}
\begin{split}
  & \hat{f}_{\mu^1,\dots,\mu^l}(t,\nu)
  = (t^{-1/2}-t^{1/2})^{l-2} \cdot
  \Big( \frac {t^{k/2}-t^{-k/2}} {t^{1/2}-t^{-1/2}} \Big)^{-2} \cdot
  \nu^{k(r-1)n/2} \cdot
  \\ & \quad \quad \quad
  \sum_{\lambda^1 \vdash n_1,\dots,\lambda^l \vdash n_l}
  \prod_{i=1}^l S_{\lambda^i,\mu^i}
  \Big( t^{(k-1)/2},t^{(k-3)/2},\dots,t^{-(k-1)/2} \Big) \cdot
  \\ & \quad \quad \quad
  \sum_{\lambda \vdash rn}
  g_{\lambda^1,\dots,\lambda^l}^\lambda(t^k) \cdot
  s_{\lambda;\nu^{-1/2}}
  \Big( t^{(k-1)/2},t^{(k-3)/2},\dots,t^{-(k-1)/2} \Big)
\end{split}
\end{equation}
for partitions $\mu^i \vdash n_i$, $i=1,\dots,l$ and
$n=n_1+\cdots+n_l$.

Let us take $T(2,k)$, $k \not\equiv 0 \pmod 2$, as an example to 
illustrate our verification of the Labastida-Mari\~no-Vafa conjecture
for torus links. In this case,
we have
\begin{equation}
\begin{split}
  & \hat{f}_{(2)}(t,\nu)
  = - \nu^k (\nu^{\frac12}-\nu^{-\frac12})^2
  (\nu+\nu^{-1}-t-t^{-1}) \cdot
  \\ & \quad \quad \quad \quad \quad \quad
  \frac {
    (t^{\frac{k+1}{2}}-t^{-\frac{k+1}{2}})
    (t^{\frac{k}{2}}-t^{-\frac{k}{2}})
    (t^{\frac{k-1}{2}}-t^{-\frac{k-1}{2}})^2
  } {
    (t^{\frac{3}{2}}-t^{-\frac{3}{2}})
    (t-t^{-1})^3
    (t^{\frac12}-t^{-\frac12})
  }, \\
  & \hat{f}_{(1,1)}(t,\nu)
  = - \nu^k (\nu^{\frac12}-\nu^{-\frac12})^2
  (\nu+\nu^{-1}-t-t^{-1}) \cdot
  \\ & \quad \quad \quad \quad \quad \quad
  \frac {
    (t^{\frac{k+1}{2}}-t^{-\frac{k+1}{2}})^2
    (t^{\frac{k}{2}}-t^{-\frac{k}{2}})
    (t^{\frac{k-1}{2}}-t^{-\frac{k-1}{2}})
  } {
    (t^{\frac{3}{2}}-t^{-\frac{3}{2}})
    (t-t^{-1})^3
    (t^{\frac12}-t^{-\frac12})
  }.
\end{split}
\end{equation}
Since $k$ is odd, both functions agree to the
Labastida-Mari\~no-Vafa conjecture. Following this way we can verify
the Labastida-Mari\~no-Vafa conjecture for all the torus links
in Examples 6.3, 6.4, and 6.5.

\end{document}